\def\versiondate{30 July 1999}
\input math.macros
\input Ref.macros

\checkdefinedreferencetrue
\continuousfigurenumberingtrue
\theoremcountingtrue
\sectionnumberstrue
\forwardreferencetrue
\initialeqmacro


\def\C{{\cal C}}
\def\P{{\bf P}}
\def\fmsf{{\ss FMSF }}

\def\Aut{{\rm Aut}}


\ifproofmode \relax \else{Version of \versiondate}\fi
\vglue20pt

\title{Nonamenable Products are not Treeable}
\author{Robin Pemantle and Yuval Peres}

\abstract{Let $X$ and $Y$ be infinite graphs, such that
   the automorphism group of $X$ is nonamenable,
and the automorphism group of $Y$ has an infinite orbit.
We prove that there is no automorphism-invariant measure 
          on the set of spanning trees in the direct product 
          $X \times Y$. This implies that
          the minimal spanning forest corresponding to i.i.d.\
          edge-weights in such a product, has infinitely many connected
          components almost surely.}

\bottom{Primary 60B99. 
Secondary 
60D05, 
20F32   
.}
{Cayley graphs, amenability, spanning trees.}
{Research partially supported NSF grant DMS-9803597.}

\bsection {Introduction}{s.intro}

There are several natural ways to pick a random spanning tree
in a finite graph,  notably the {\it uniform spanning tree}
and the {\it minimal spanning tree} (for random edge weights).
The limits of these models on infinite graphs
sometimes give spanning {\it forests}, rather than trees.
In this note, we present a large class of graphs where
one cannot pick a random subtree in an automorphism-invariant 
manner. 

\noindent{\bf Definition.}\enspace
 A subtree $\Gamma=(V_\Gamma,E_\Gamma)$ in a graph $X=(V_X,E_X)$ is an acyclic 
connected graph, with $V_\Gamma \subset V_X$
and $E_\Gamma \subset E_X$.  If $V_\Gamma=V_X$, 
then $\Gamma$ is called a {\it spanning\/} tree of $X$.
By identifying a subtree with the indicator function of its edges,
we may view the ensemble of spanning trees in $X$ as a Borel set
in the compact metrizable space $\{0,1\}^{E_X}$.

\procl t.span
   Let $X$ and $Y$ be infinite, locally finite graphs. 
Suppose that $G \subset \Aut(X)$ 
is a closed nonamenable subgroup of $\Aut(X)$,
and $H \subset \Aut(Y)$ has an infinite orbit.
    Then there is no $G \times H$-invariant
    probability measure on the set of
       spanning trees of the direct product graph $X \times Y$. 
\endprocl

Given two graphs $X=(V_X,E_X)$ and $Y=(V_Y,E_Y)$, 
the {\bf direct product graph} $X \times Y$ has vertex set $V_X \times V_Y$;
the vertices $(x_1,y_1)$ and $(x_2,y_2)$
in $V_X \times V_Y$ are taken to be adjacent in $X \times Y$ iff either 
$x_1=x_2$ and $[y_1,y_2] \in E_Y$, or $y_1=y_2$ and 
$[x_1,x_2] \in E_X$.

In the next section we provide  some probabilistic motivation for \ref t.span/,
by describing its application to {\it minimal spanning trees\/} 
(corresponding to i.i.d.\ weights) in finite graphs, and their
limits, {\it minimal spanning forests\/} in infinite graphs.

The following extension of \ref t.span/ is sometimes useful.
\procl c.span
   Under the hypothesis of \ref t.span/, there is no $G \times H$-invariant
    probability measure on the set of
       subtrees of the direct product graph $X \times Y$. 
\endprocl
\proof
>From the assumption that $\Gamma$ is a random subtree in $X \times Y$
with a $G \times H$-invariant law, we will obtain a contradiction.
 Let ${\cal S}(\Gamma,n)$
denote the set of vertices in $X \times Y$ at graphical distance $n$ from $\Gamma$.
For every $n \ge 1$ and every $v \in {\cal S}(\Gamma,n)$,
choose uniformly at random a neighbor $v'$ of $v$ such that  
$v' \in {\cal S}(\Gamma,n-1)$, and add the edge $[v,v']$ to $\Gamma$.
This yields a random spanning tree in $X \times Y$
with a $G \times H$-invariant law, and contradicts \ref t.span/. 

\Qed

Next, we present a variant of \ref t.span/ that does not require any apriori graph structure.

\noindent{\bf Definitions.}\enspace
 {\bf (i)} A spanning tree in a countable set $\Upsilon$ is an acyclic 
connected graph
with vertex set $\Upsilon$; it need not be locally finite.
We may view the ensemble of spanning trees in $\Upsilon$ as a Borel set
in the compact metrizable space $\{0,1\}^{\Upsilon\times \Upsilon}$.

\noindent{\bf (ii)} A countable group $G$ is {\bf treeable} 
if there exists a Borel
probability measure on the set of spanning trees in $G$, which is
invariant under the action of $G$ by right multiplication.

\smallskip
Obviously, a free group is treeable. The work of Ornstein and Weiss (1987) 
 implies that any discrete amenable group is treeable;
 see Theorem 5.3 of Benjamini, Lyons, Peres and Schramm (1999a), 
denoted BLPS (1999a)
below, for an alternative proof.
On the other hand,
 results of Adams and Spatzier (1990) imply that 
 groups with Kazhdan's property T are not treeable.
The following theorem exhibits simpler examples of nontreeable groups.
\procl t.span-group
   Let $G$ and $H$ be countably infinite groups. If the direct product
$G \times H$ is treeable, then $G$ and $H$ are amenable. 
\endprocl
In particular, the direct product of a free group 
(on two or more generators) 
with any infinite group is nontreeable. 
It is an intriguing unsolved problem to find a geometric 
characterization for treeable groups.

\ref t.span-group/ seems close in spirit to results of Adams (1988, 1994);
however, it appears that there is no direct implication from one 
result to the other
(S. Adams, R. Lyons and B. Weiss, private communications). 
The referee has suggested that perhaps the {\it techniques}
of Adams (1988) could be used to give an alternative proof
of our main results.

After seeing an earlier version of
the present paper, R. Lyons has pointed out that 
\ref t.span-group/ (but not \ref t.span/)
could also be inferred from recent work of Gaboriau (1998)
on the ``cost'' of equivalence relations.


In the next section we discuss minimal spanning trees and forests.
\ref s.back/ contains background material on amenability.
Finally, in \ref s.free/, we establish a
common generalization of 
\ref t.span-group/ and \ref t.span/.

\bsection {Minimal spanning trees and forests}{s.mst}
Suppose that the edges of a connected finite graph
  $(V,E)$ are labeled with i.i.d.\ random variables $\{U_e\}_{e \in E}$,
that are uniform in $[0,1]$. The corresponding {\bf minimal spanning tree} $\Gamma$
(the spanning tree that minimizes the sum of labels on its edges)
consists of  edges $e$ such that 
there is no path between the endpoints of $e$ where all edges 
on the path have labels lower than $U_e$.

If we now consider a sequence of connected finite graphs $X_n$ that
exhaust an infinite, locally finite connected graph $X$, 
then the corresponding minimal
spanning trees $\Gamma_n$ converge almost surely to a random spanning subgraph
of $X$ called the {\bf Free minimal spanning forest} (\fmsf).
See Alexander (1995) and BLPS (1999c) for more information
on minimal spanning forests, and their connections to percolation.

The \fmsf $F_X$ on $X$ may be constructed directly, as follows:
Label the edges of $X$ by i.i.d. random variables
$\{U_e\}$, uniform in $[0,1]$. Then
 remove any edge $e$ that has the highest label
in some finite cycle. In other words, $e$ is retained in $F_X$ iff
there is no path between the endpoints of $e$ consisting of edges with
labels lower than $U_e$.
It is clear that $F_X$ has no cycles a.s. Moreover, all connected
components of $F_X$ are infinite a.s., since for any finite set
of vertices $K$, the edge $e$ that has the lowest label
among the edges connecting $K$ to its complement, must be in $F_X$. 

Newman and Stein (1996) conjectured that the \fmsf 
in $\Z^d$ is disconnected if $d$ is large.
This is still open, but \ref t.span/ implies the following.

\procl c.msf
     Let $X$ be a nonamenable connected graph $X$ that has a 
 quasi-transitive unimodular automorphism group $G$, and
 suppose that the graph $Y$ has an automorphism group $H$ 
   with an infinite orbit. Then the \fmsf on $X \times Y$ 
 has infinitely many components a.s.
\endprocl

\proof
Denote by $N(F)$ the number of components of the \fmsf $F$ on $X \times Y$.
Since $N(F)$ is invariant under the ergodic action of $G \times H$,
it is a.s.\ constant. \ref t.span/ ensures this constant is 
greater than $1$. Moreover, if $N(F)<\infty$, then choosing uniformly one of the $N(F)$
components of $F$ would yield a $G \times H$-invariant measure on subtrees
of $X \times Y$, contradicting  \ref c.span/. Therefore, $N(F)=\infty$ a.s.
\Qed
We note that the
{\it uniform\/} spanning forest in a nonamenable
 product $X \times Y$ as in \ref t.span/ is known
to have infinitely many components, see BLPS (1999b).

The next corollary concerns minimal spanning trees in certain 
finite graphs. Nevertheless, we do not know any finitistic proof.

\procl c.mst
Let $T_d(n)$ be a finite tree with root $\rho$, where every vertex has degree
$d \ge 2$, except the vertices at the maximal distance $n$ from the root,
that have degree $1$. 
 Let $\Gamma_n$ be the minimal spanning tree in the product
graph $T_d(n) \times T_b(n)$,  determined by i.i.d.\ uniform labels 
on the edges.
Denote by $L_n$ the sum of the distances in $\Gamma_n$ from $(\rho,\rho)$
to its neighbors in $T_d(n) \times T_b(n)$. If $d \ge 3$,
then the random variables $L_n$ are not tight, {\it i.e.},
$\sup_{M>1} \inf_{n} \P[L_n \le M] <1$.
\endprocl

\proof 
Denote the $d$-regular infinite tree by $T_d$, and consider
the finite trees $T_d(n)$ as embedded in $T_d$, with the roots identified.
Label the edges of  $T_d \times T_b$ by uniform
variables $\{U_e\}_{e \in E}$. 
The resulting minimal spanning trees on $T_d(n) \times T_b(n)$
converge almost surely to the \fmsf $F$ in $T_d \times T_b$.

If the variables $L_n$ were tight, then passage to the limit 
would imply that the distance in $F$ between any two vertices
that are adjacent in $T_d \times T_b$ 
is finite a.s., so $F$ is connected, contradicting \ref c.msf/.
\Qed

\bsection{Background on amenability}{s.back} 
Let $G$ be a locally compact group, endowed
 with a left-invariant Haar measure.
A linear functional $M$ on $L^\infty(G)$ is called a 
{\bf mean} if it maps the constant
function $\I{}$ to 1 and $M(f) \ge 0$ for $f \ge 0$.
 If
$f \in L^\infty(G)$ and $g \in G$, we write $L_g f(h):= f(g h)$. We
call a mean $M$ {\bf invariant} if $M(L_g f) = M(f)$ for all $f
\in L^\infty(G)$ and $g \in G$. Finally, we say that $G$ is {\bf amenable} if
there is an invariant mean on $L^\infty(G)$. 
See Paterson (1988) for properties and characterizations of amenable groups,
and the papers Soardi and Woess (1990), Salvatori (1992) and
BLPS (1999a) for the relation between isoperimetric
inequalities in a graph $X$ and nonamenability of 
$\Aut(X)$.

We need the following variant of a method due to Adams and Lyons (1990).
\procl l.adly
 Let $G$ be a locally compact group, that acts
transitively on the countable set $X$, so that for every $x \in X$,
the {\bf stabilizer} 
$
S_x:=\{g \in G \st g(x)=x\} 
$ is compact.
We are given that $G$ also acts by measure-preserving
maps on a probability space $(\Omega,{\cal F}, \P)$.
Suppose that for each $n \ge 1$
 and $\omega \in \Omega$, an
equivalence relation
$R_n=R_n(\omega) \subset  X \times X$ is given,
so that for any $x_1,x_2 \in X$, the set
$\{\omega  \st (x_1,x_2) \in R_n(\omega)\}$ is in ${\cal F}$.
We assume that 
$$
\forall g \in G, \qquad (x_1,x_2) \in R_n(\omega) \; \hbox{ \rm iff } \; 
(g x_1,g x_2) \in R_n(g \omega) \, . \label e.rinv
$$
If for every $x \in X$ and  $n \ge 1$, and for almost every 
$\omega \in \Omega$, the 
equivalence class $\C_n(x):=\{z \in X \st  (x,z) \in R_n\}$
is finite, 
and $\lim_n \P[(x,z)\in R_n]=1$ for every $x,z \in X$,
then $G$ is amenable.
\endprocl
\proof
The hypothesis that $G$ preserves $\P$, together with \ref e.rinv/,
gives
$$
\P[(gx,gz) \in R_n] = \P[(x,z) \in R_n]
\label e.rinv2
$$
for all $x,z \in X$ and $g \in G$.

A linear functional $M$ on $\ell^\infty(X)$ is called a 
{\it mean} if it maps the constant 
function $\I{}$ to 1 and $M(f) \ge 0$ for $f \ge 0$.
 For $f \in \ell^\infty(X)$ and $g \in G$, we write $L_g f(x):= f(g x)$. 
A mean $M$ on $X$ is {\bf $G$-invariant} if $M(L_g f) = M(f)$ for all 
$f \in \ell^\infty(X)$.

Fix $o \in X$. For each $n$, define a mean $M_n$ on $X$ by
$M_n(f):= \E[\sum_{x \in \C_n(o)} f(x) / |\C_n(o)|]$.
By \ref e.rinv2/,
$$
M_n(L_g f)=\E[\sum_{g^{-1}z \in \C_n(o)} f(z) / |\C_n(o)|]=
\E[\sum_{z \in \C_n(g o)} f(z) / |\C_n(go)|] \,.
$$
Therefore 
$$
M_n(L_g f)-M_n(f) \le 2\|f\| \P[(o,go) \notin R_n] \to 0 
\hbox{ as } n \to \infty \,.
$$
Consequently, any weak$^*$ limit point $M_*$ of $\{M_n\}$
is a $G$-invariant mean on $X$.
Next, given $f \in L^\infty(G)$, define $\overline f \in
\ell^\infty(X)$ by 
$
\overline f (x) := \mu(S_o)^{-1}  \int_{ho = x} f(h)\,d\mu(h) 
$,
where $\mu$ is Haar measure. It is easy to check that
$\overline{L_g f} = L_g \overline f$, so that an invariant mean $M$ on $G$
may be defined by $M(f):= M_*(\overline f)$.
\Qed

\bsection{Spanning trees not confined to graphs}{s.free}
The following theorem extends \ref t.span/ and \ref
t.span-group/. 
\procl t.span2
   Let $G$ be a locally compact group, that acts 
on the countable set $X$, so that
 for each $x \in X$, the stabilizer $S_x \subset G$ is compact.
       Let $H$ be a group that acts on the countable set $Y$
so that 
$$
\hbox{for any finite subset $Y'$ of $Y$, there exists $h \in H$
satisfying $h(Y') \cap Y' = \emptyset$.}
\label e.slug
$$
 If there exists a $G \times H$-invariant probability measure on 
      the ensemble of  spanning trees in $X \times Y$, then $G$ is amenable.
\endprocl

\noindent{\bf Remark.} 
Let $H$ be  a group of permutations of the countable set $Y$.
As noted by B. Weiss and the referee,
the hypothesis \ref e.slug/ on $H$ holds if
and only if all $H$-orbits in $Y$ are infinite. See
Newman (1976).

The next lemma is obvious if $G$ is a group of graph automorphisms.
\procl l.balls
Suppose that the locally compact group $G$  acts 
on the countable set $X$ with compact stabilizers. Then there exist
sets $\{B(x,\ell) \st x \in X, \ell \ge 1\}$ in $X$, such that for any $x \in X$
and $g_1 \in G$, we have
$$ 
g_1 B(x,\ell)=B(g_1 x,\ell) \; \hbox{ \rm and } \,  \cup_{\ell \ge 1}  B(x,\ell) =X \, .
\label e.bprop
$$
Moreover, for  any $x_1 \in X$ and $\ell \ge 1$, the set
$\{x_2 \in X \st B(x_1,\ell) \cap B(x_2,\ell) \neq \emptyset\}$ is finite.
\endprocl
\proof
Fix $x_0 \in X$, and finite sets $\{X_\ell\}_{\ell \ge 1}$
 such that $\cup_{\ell \ge 1}  X_\ell =X$. Define
$$
B(x, \ell):=\{g z \st z \in X_\ell, \, g x_0=x\} \,.
$$
This is a finite set, since  the sets $\{g   \st  g x_0=x, gz=w\}$ 
with $z \in X_\ell$ and $w \in B(x, \ell)$, form an open cover of the compact set
 $\{g   \st  g x_0=x\}$.
 The  properties in \ref e.bprop/ are immediate.
Finally, for any $w \in X$,
$$
\{x_2 \in X \st w \in B(x_2,\ell) \}= \bigcup_{z \in X_\ell}
\{g x_0 \st gz=w\}
$$
is a finite set. Taking the union of these sets over $w \in B(x_1,\ell)$
completes the proof of the lemma.
\Qed

\proofof t.span2
Let $\Upsilon:=X \times Y$ and denote by 
$\Omega \subset \{0,1\}^{\Upsilon \times \Upsilon}$
the set of indicator functions of spanning trees in $\Upsilon$.
In particular, $\omega(u,v)=\omega(v,u)$ for any  $\omega \in \Omega$
and $u,v \in \Upsilon$.
Our hypothesis is that that there is a $G \times H$-invariant
probability measure $\P$ on $\Omega$.

For $\omega \in \Omega$, denote by $\Gamma(\omega)$ the corresponding
spanning tree, {\it i.e.}, the set of unordered
pairs
$\{u,v\} \subset \Upsilon$ such that $\omega(u,v)=1$.
Let $\{Y_n\}_{n \ge 0}$ be an increasing sequence of finite sets, such that
$Y_0=\{y_0\}$ and $\cup_n Y_n=Y$.
For each $n$, pick $h_n \in H$ such that $h_n(Y_n) \cap Y_n =\emptyset$,
and denote $y_n=h_n(y_0)$. 
Recall the sets $B(x,\ell)$ from \ref l.balls/;
the parameter $\ell=\ell(n)$ will be specified below.

Fix an orbit $X_o$ of $G$ in $X$.
For $\omega \in \Omega$ and $n \ge 1$,
define the equivalence relation $R_n=R_n(\omega)$ on $X_o$,
by letting $(x_1,x_2) \in R_n$ iff
\item{(i)}
 the path in $\Gamma(\omega)$ from $(x_1,y_0)$ to $(x_2,y_0)$ 
is contained in $X \times Y_n$, and
\item{(ii)}
 the path in $\Gamma(\omega)$ from $(x_1,y_n)$ to $(x_2,y_n)$ 
is contained in $X \times h_n(Y_n)$;
\item{(iii)} the path in $\Gamma(\omega)$ from $(x_1,y_0)$ to $(x_1,y_n)$ 
is contained in $B(x_1, \ell)  \times Y$;
\item{(iv)}the path in $\Gamma(\omega)$ from $(x_2,y_0)$ to $(x_2,y_n)$ 
is contained in $B(x_2, \ell)  \times Y$.

For fixed $x_1,x_2$ in $X_o$, the events in (i) and (ii) above have the same 
probability, which tends to $1$ as $n \to \infty$. 
By choosing $\ell=\ell(n)$ large enough, we can ensure that the events
in (iii) and (iv) have probability at least $1-1/n$.
for $x_1,x_2$ in $X_o$.
Therefore
$\lim_n \P\Big( (x_1,x_2) \in  R_n \Big) =1$ for any $x_1,x_2$ in $X_o$.
The invariance relation \ref e.rinv/ is easily checked,
so in order to apply \ref l.adly/, we just need to verify that
for every $x_1 \in X_o$, the equivalence class $\C_n(x_1)$ is a.s.\ finite.

Suppose that $(x_1,x_2) \in  R_n$ and consider the ``cycle''
obtained by concatenating the following paths in $\Gamma(\omega)$:
$$ 
(x_1,y_0) \to (x_1,y_n) \to (x_2,y_n) \to (x_2,y_0) \to (x_1,y_0) \,.
\label e.con
$$
Since $\Gamma(\omega)$ is a tree, every edge that is traversed in \ref e.con/, 
must be traversed an even number of times. Parts (i) and (ii) in the definition of $R_n$
imply that the first edge in the  $\Gamma(\omega)$-path 
$(x_1,y_0) \to (x_1,y_n)$
that exits $X \times Y_n$, must also occur in the 
$\Gamma(\omega)$-path  $(x_2,y_n) \to (x_2,y_0)$.
Therefore, by parts (iii) and (iv) of that definition, for $\ell=\ell(n)$ we have
$ B(x_1,\ell) \cap B(x_2,\ell) \neq \emptyset$.
Thus \ref l.balls/ ensures that
$\C_n(x_1)$ is a.s.\ finite.
We have verified all the hypotheses of \ref l.adly/
(with $X_o$ in place of $X$), so $G$ is amenable.

\Qed

\medbreak
\noindent {\bf Acknowledgement.}\enspace 
We are grateful to 
Itai Benjamini, Russ Lyons and Oded Schramm for numerous 
discussions that led to this work. Russ Lyons suggested that a unimodularity
assumption in a previous version of this note could be dispensed with.
We also thank the referee, Scot Adams and Benjy Weiss for useful remarks.

\beginreferences


 Adams, S. (1988) Indecomposability of treed equivalence relations. 
{\it Israel J. Math.} {\bf 64}, 362--380.

 Adams, S. (1994) Indecomposability of  equivalence relations
generated by word-hyperbolic groups.
{\it Topology} {\bf 33}, 785--798.

Adams, S. \and Lyons, R. (1991) Amenability, Kazhdan's property and
percolation for trees, groups and equivalence relations,
{\it Israel J. Math.} {\bf 75}, 341--370.

Adams, S. \and Spatzier, R.~J. (1990) Kazhdan groups, cocycles and trees. 
{\it Amer. J. Math.} 112 (1990), 271--287.

Alexander, K. S. (1995) Percolation and minimal spanning
forests in infinite graphs, {\it Ann. Probab.} {\bf 23}, 87--104.

Benjamini, I., Lyons, R., Peres, Y., \and Schramm, O. (1999a)
Group-invariant percolation on graphs, 
{\it Geom.\ Funct.\ Anal.} {\bf 9}, 29--66. 

Benjamini, I., Lyons, R., Peres, Y. \and Schramm, O. (1999b)
Uniform spanning forests, {\it Ann. Probab.}, to appear.

Benjamini, I., Lyons, R., Peres, Y. \and Schramm, O. (1999c)
Minimal spanning forests, {\it in preparation}.

Gaboriau, D. (1998) Co\^ut des relations d'\'equivalence et des groupes.
{\it Preprint.}

H\"aggstr\"om, O. (1997) Infinite clusters in dependent automorphism
invariant percolation on trees, {\it Ann. Probab.} {\bf 25}, 1423--1436.





Newman, C. M. \and Stein, D. L. (1996) Ground-state structure
in a highly disordered spin-glass model, {\it J. Statist. Phys.} {\bf 82},
1113--1132.

Newman, P. M. (1976) The structure of finitary permutation groups,
{\it Arch.\ Math.} (Basel) {\bf 27}, 3--17.

Ornstein, D. S. \and Weiss, B. (1987) Entropy and isomorphism theorems
for actions of amenable groups, {\it J. d'Analyse\/} {\bf 48}, 1--141.

Paterson, A. L. T. (1988) {\it Amenability}. American Mathematical Soc.,
Providence.

Salvatori, M. (1992) On the norms of group-invariant transition operators
on graphs, {\it J. Theor. Probab.} {\bf 5}, 563--576.


Soardi, P. M. \and Woess, W. (1990) Amenability,
unimodularity, and the spectral radius of random walks on infinite  
graphs, {\it Math. Z.} {\bf 205}, 471--486. 

Trofimov, V. I. (1985) Groups of automorphisms of graphs as topological
groups, {\it Math. Notes} {\bf 38}, 717--720.


\endreferences

\filbreak
\begingroup
\eightpoint\sc
\parindent=0pt\baselineskip=10pt
\def\email#1{\par\qquad {\tt #1} \smallskip}
\def\emailwww#1#2{\par\qquad {\tt #1}\par\qquad {\tt #2}\smallskip}

Department of Mathematics, University of Wisconsin, Madison, WI 53706 
\email{pemantle@math.wisc.edu}

\medskip

Institute of Mathematics,
The Hebrew University, 
Givat Ram, Jerusalem 91904,
Israel 

and Department of Statistics, University of California, Berkeley, CA
\emailwww{peres@math.huji.ac.il}
{http://www.ma.huji.ac.il/\~{}peres/}

\bye